\documentclass[a4paper,12pt]{amsart}

\usepackage{amsmath, amsfonts, amsthm, amsfonts}

\def\]{\textup{\mbox{]\hspace{-.15em}]}}}
\def\[{\textup{\mbox{[\hspace{-.15em}[}}}

\def\got{\mathfrak}

\newcommand{\Ind}{\mathrm{Ind}}
\newcommand{\W}{\mathrm{W}}

\newcommand{\Dpst}{{\rm D_{pst}}}

\newcommand{\Fur}{F_v^{ur}}
\newcommand{\Ker}{\mathrm{Ker}}
\newcommand{\Q}{\mathbb{Q}}
\newcommand{\R}{\mathbb{R}}
\newcommand{\Qpb}{\overline{\Q_p}}

\newcommand{\U}{\mathrm{U}}
\newcommand{\UU}{\mathcal{U}}
\newcommand{\GL}{\mathrm{GL}}
\newcommand{\GSp}{\mathrm{GSp}}
\newcommand{\AAA}{\mathbb{A}}

\newcommand{\ps}{\par \medskip }
\newcommand{\pn}{\par \noindent}
\newcommand{\Z}{\mathbb{Z}}
\newcommand{\C}{\mathbb{C}}

\newcommand{\OO}{\mathcal{O}}

\newcommand{\Gal}{\mathrm{Gal}}

\newcommand{\pf}{ {\it Proof.\,\,  } }
\newcommand{\epf}{$\square$ \par }

\newcommand{\G}{\mathrm{G}}
\newcommand{\rond}{\circ}

\newtheorem{thm}{Theorem}
\newtheorem{prop}{Proposition}
\newtheorem{cor}{Corollary}
\newtheorem{lemme}{Lemma}

\title{On number fields with given ramification}

\author{Ga\"etan Chenevier}

\begin{document}

\maketitle

\section*{Introduction}\label{intro}

        Let\footnote{December 2004.} $E$ be a number field and $S$ a nonempty set of places of
$E$. We denote by $E_S$ a maximal algebraic extension of $E$ unramified outside $S$.
Let us fix $u \in S$ and an $E$-embedding $\varphi$ of $E_S$ in an algebraic closure $\overline{E_u}$ of
$E_u$. In this paper, we are interested in the following property: \ps

\begin{center}$E_S$ is dense in $\overline{E_u}$, \end{center} 

\vspace{-0,7cm} \begin{flushleft} ($P_{E,S,u}$) \end{flushleft} 

\bigskip

\noindent where the identification of $E_S$ with $\varphi(E_S)$ is understood. \noindent It is easy to see that
$P_{E,S,u}$ is independent of the choice of $\varphi$, and equivalent to each of the following
properties:\ps

(i) the map $\Gal(\overline{E_u}/E_u) \longrightarrow \Gal(E_S/E)$
induced by $\varphi$ is injective,\ps
(ii) the $E_u$-vector space generated by $E_S$ is $\overline{E_u}$, \ps
(iii) for all finite extension $K/E_u$, there exists a number
field $E'/E$ unramified outside $S$ and a place $u'|u$ such that $K$
has a continuous embedding into $E'_{u'}$. \ps

A trivial remark is that if $u \in S' \subset S$, then $P_{E,S',u}$ implies $P_{E,S,u}$. Moreover, a simple argument using Krasner's lemma and a weak approximation theorem shows that
a stronger form of $(P_{E,S,u})$ is true\footnote{Let $M$ be a local field, $d \geq 1$ 
an integer, and $Q \in M[T]_d$ a separable polynomial of degree $d$. By continuity 
of roots and Krasner's lemma, if $R\in M[T]_d$ is 
sufficiently close to $Q$, then $M[T]/(R) \simeq M[T]/(Q)$ as $M$-algebra. Now, let $S'$ be a 
finite set of places of a number field $E$, and let us fix for each $x \in S'$ an
\'etale $E_x$-algebra $A_x$, each of same degree $d$. Then the argument above, the primitive element theorem,
and the weak approximation for the affine line over $E$, show that there exists an
\'etale $E$-algebra 
$A$ of degree $d$, with $A \otimes_E E_x \simeq A_x$ for all $x \in S'$. It implies the claim 
if we take $S'$ to be the set of places $x$ of $E$ such that $x=u$ or $x \notin S$, $A_u:=K$, 
and (for example...) $A_x=E_x^d$ if $x\in S'\backslash \{u\}$. Of course, this simple approach 
is inefficient when $S'$ is not finite.} if $S$ contains almost all the places of $E$. 
However, Hermite's theorem asserts that $\Q_{\{\infty\}}=\Q$, hence
$P_{\Q,\{\infty\},\infty}$ is not satisfied. For these reasons and others 
related to arithmetic geometry which will become clear later, we shall focus our 
attention on the case where $S$ is a {\it finite set} containing all the 
archimedian places of $E$ and all the finite places dividing a given prime
number $l$. As $P_{E,S,u}$ is then obviously true for each archimedian place $u$
(e.g. because of cyclotomic fields), let us assume also that $u$ is finite. At least in this setting, the question
of deciding if $P_{E,S,u}$ is true seems to be in the folklore of algebraic 
number theory\footnote{In particular, it had been asked by Ralph Greenberg after a related question raised by James Milne concerning the procardinal of $\Gal(E_S/E)$ (see \cite[Chap. I \S 4]{Mi} and our corollary \ref{cor1}).}, and as far as the author knows, there was no 
example before this paper of such a triple $(E,S,u)$ where the answer was known. Note that although local 
Galois groups are pro-solvable, it does not seem possible\footnote{The author
confesses that he did not manage to prove that the "the maximal
pro-solvable extension of $E$ inside $E_S$ is not dense in
$\overline{E}_u$", although it seems reasonable.} to deduce property 
$P_{E,S,u}$ by induction from class field theory (or by Grunwald-Wang theorem, 
\cite{AT}[p.105]), the obstructions given by units forcing us to enlarge $S$ 
at each step (see also prop. \ref{general} \S \ref{constructionII}). Let us mention also that a weak version of $P_{E,S,u}$ is known by the work of Kuz'min and V.G. Mukhamedov (\cite{MU}, 
\cite{NSW} chapter X, \S 6, thm. 10.6.4 and the last exercise) concerning
$p$-extensions of CM fields. For example, the following result is proved in \cite{MU}. Let 
$p$ be a prime number and let $E$ be a CM field with maximal totally real subfield $F$ such that each prime $v$ of $F$ dividing $p$ splits in
$E$. If $v$ is such a place, then the canonical maps $$\Gal(\overline{F_v}/F_v)_p \longrightarrow
\Gal(\overline{E}_{\{p,\infty\}}/E)_p,$$ \noindent are injective. In this statement, $H_p$ denotes the maximal
pro-$p$-quotient of the profinite group $H$. However, there seems to be no way to deduce
property $P_{E,\{\infty,p\},v}$ from these ones. We come now to our results: \ps

\begin{thm} \label{T1} Assume that $E$ is a $CM$ field and that $u$ is split above a finite place $v$ of
the maximal totally real subfield $F$ of $E$. If $l$ is a rational prime number which is prime to $v$, 
and if $S$ is the set of places of $E$ dividing $lv$, then $P_{E,S,u}$ holds. 
\end{thm}

\begin{cor} \label{cor1}If $E$ and $S$ are as in theorem \ref{T1}, then any integer $n\geq 1$ divides
the pro-cardinal of $\Gal(E_S/E)$.\end{cor}

Indeed, this last property is a general consequence of $P_{E,S,u}$ when $u$
is finite, as $\Gal(\overline{E_u}/E_u)$ has then a continuous
surjective homomorphism to $\widehat{\Z}$. As Milne pointed out to us, this corollary answers in some cases a question raised in \cite[Chap. I \S 4]{Mi} (concerning the set $P$ defined there). Note that the hypothesis of theorem \ref{T1} are satisfied for instance if $E$ is a 
quadratic imaginary field split at a prime number $v=p$, and $l$ is a prime $\neq p$. It has the following
consequence for $E=\Q$:

\begin{cor} \label{cor2} Let $p$ be a prime number, $N$ an integer such that $-N$ is
the discriminant of an imaginary quadratic field in which $p$ splits, and
let $S$ be the set of primes dividing $Np\infty$. Then $P_{\Q,S,p}$ holds. \end{cor}

\noindent Explicit examples are given by $(p,N) \in \{ (2,7), (3,2^3), (3,11), (5,2^2), (5,11), (7,3), \dots
\}$. In fact, as we will see later, we certainly conjecture that if $p$ and $l$ are distinct
primes, then $P_{\Q,\{\infty,p,l\},p}$ holds:

\begin{thm}\label{T2} Let $p$ and $l$ be any distinct prime numbers, and assume 
hypothesis $1$ and $2$ (see \S \ref{conj}), then $P_{\Q,\{\infty,p,l\},p}$ holds.
\end{thm}

Although it is tempting to conjecture that $P_{\Q,\{\infty,l\},l}$ holds, we are unfortunately less definite
in this case (see \S \ref{l=p}).

\medskip

Let us describe a bit the strategy of the proofs. First, we have to construct a lot of number fields
unramified outside a given set of places $S$. By a well known result of
Grothendieck, the number fields attached to the $l$-adic \'etale cohomology of
a proper smooth scheme $X$ over $E$ satisfy this property if $X$ has good
reduction outside $S$, and if $S$ contains the primes dividing $l$. Although it
might be very difficult in general to find such an $X$ ($S$ being given),
well chosen Shimura varieties give some interesting examples.
%\footnote{For instance, $X_1(l^n)$ has good reduction outside $l$, and its genus goes to $\infty$ with $n$. However, it only constructs number fields of $\GL_2$-type in some sense, hence they cannot exhaust $\overline{\Q_l}$.}
Even better, their $l$-adic cohomology is completely described, at least conjecturally, by the Langlands conjectures, in terms of cohomological automorphic forms. By the work of many authors, culminating to Harris and Taylor's proof of
local Langlands conjecture for $\GL_n$, a big part of these conjectures is known for the
so-called "simple" Shimura varieties, which are attached to some unitary groups. This combined with some other 
tricks allowed Harris and Taylor to attach an $l$-adic representation to a wide class of cuspidal automorphic forms for 
$\GL_n(\AAA_E)$ (\cite[thm. C]{HT}), compatible with the local Langlands correspondence at all finite places not dividing $l$. Thus, we first show that we can construct cuspidal automorphic forms
for $\GL_n(\AAA_E)$ satisfying Harris-Taylor's conditions, and which are unramified outside $S$ and of given "ramification type" at $u$. More precisely, we first construct some automorphic forms on well-chosen unitary groups and we apply quadratic base-change to them (\S \ref{unit}, \S \ref{auto}). By Harris and Taylor's
result, their associated $l$-adic Galois representations give us number fields with the
required local properties and some control at $u$ (\S \ref{HT}). The local Langlands correspondence shows then that 
we can produce in this way a lot of such number fields, using unitary groups {\it of all ranks}. Two little subtleties arise at this point. First, we have few control, of course, on the Weil numbers in the above constructions,
i.e. on the {\it unramified part} of the completion at $u$ of the number fields
constructed. The second one is that the automorphic representations we consider
must fulfill some conditions so that we may apply to them the results of Harris and Taylor and the known cases of quadratic base change. However, a simple trick (\S
\ref{trick}) allows us to show that we produced sufficiently many number
fields to prove theorem \ref{T1}. \par
As the sketch above shows, a question closely related to our initial aim is to ask if there
exists pure motives over $E$ with "reduction type" prescribed at
each finite place (and say generic Hodge numbers). On the automorphic side, 
it leads to the purely analytic problem of constructing non trivial, discrete, algebraic 
automorphic forms for a given reductive $\Q$-group $H$, with prescribed 
properties at all finite places. Of course, it is not always possible, e.g. an
anisotropic torus has this property if and only if its real points are
compact. In \S \ref{constructionII}, we show how to construct these automorphic
forms under the following hypothesis on $H$: its center has finite arithmetic subgroups 
and $H(\R)$ has a holomorphic discrete series. This result was probably
known to specialists, but we could not find any convenient reference
in this generality. By similar arguments as in the proof of theorem \ref{T1}, and using 
these automorphic representations, we explain in the last section some consequences 
concerning property $P_{E,S,u}$ of some standard conjectures in the arithmetic theory of 
automorphic forms. In particular, we get theorem \ref{T2}. \ps

{\it Notations:} If $F$ and $G$ are two subfields of a given field, 
we denote by $F.G$ the {\it subfield} generated by $F$ and $G$. This is  
also the $F$-vector space generated by $G$ if $G$ is algebraic over $F \cap
G$. If $K$ is a field, $\overline{K}$ denotes a separable algebraic closure of $K$, and 
$\G_K:=\Gal(\overline{K}/K)$ its absolute Galois group, equipped with its Krull topology. If $\rho: \G_K \rightarrow H$ is any continuous group homomorphism to a topological group $H$, we will denote by $K(\rho):=\overline{K}^{\Ker(\rho)}$ the algebraic normal extension 
of $K$ fixed by $\Ker(\rho)$. We will say that $K(\rho)$ is the extension of $K$ {\it cut out} by $\rho$. If $K$ is a finite extension of $\Q_p$, $K^{ur} \subset \overline{K}$ denotes
its maximal unramified extension, 
$\W_K \subset \G_K$ its Weil group, and
$I_K=\Gal(\overline{K}/K^{ur}) \subset  \W_K$ its inertia group. If $F$ is a number field, 
$\AAA_F$ denotes the ad\`ele ring of $F$ and $\AAA_{F,f}$ its quotient of finite ones. If 
$p+q=n$, we denote by $\UU(p,q)$ the real unitary group
of signature $(p,q)$, and we set $\UU(n):=\UU(n,0)$.

\tableofcontents

\vspace{-3mm}
\section{Construction of some automorphic forms I} \label{constructionI}
\medskip

Let $E$ be a $CM$ field as in the statement of theorem \ref{T1}. Let us fix a
finite place $w$ of $F$ dividing the prime number $l$ (in particular $w\neq v$), and let $n \geq 1$ be an integer. \ps

\subsection{Some unitary groups} \label{unit}We want to consider unitary
groups attached to central division algebras over $E$, which are quasisplit at each finite place $\neq w,\, v$,
and compact at infinite places for convenience. The relevant "Hasse principle" is known and
due to Kneser. We refer to Clozel's paper \cite{C} \S 2 for a convenient exposition of Kottwitz' interpretation of Kneser's results
in the special case of unitary groups. \ps

\begin{lemme} \label{groupe} There exists a unitary
group $\U(n)$ over $F$ attached to $E/F$ such that for a place $x$
of $F$, $\U(n)(F_x)$ is:\ps

(a) quasisplit if $x$ is a finite place not dividing $wv$, \ps
(b) the group of units of a central division algebra over $F_x$ if $x=v$.\ps
(c) the compact real unitary group if $x$ is real.

\end{lemme}

\pf If $n$ is odd, there is no global obstruction to the existence of
such groups by \cite{C}[lemma 2.1]. Assume $n$ is even. By {\it loc.cit.}
(2.2), the global obstruction
lies in $\Z/2\Z$ and is the sum of all local ones modulo $2$. Assuming given
local groups satisfying (a), (b), and (c), we can make the global invariant
vanish
by requiring, if necessary, that $\U(n)(F_w)$ is either a non-quasisplit
unitary group or the units of a
division algebra, because such groups have local invariant $\equiv 1 \bmod
2$ by loc.cit.
(2.3). This concludes the proof. \epf

\subsection{Construction of automorphic forms} \label{auto} Let $H/F$ be the
unitary group $\U(n)$ given by lemma \ref{groupe}. The group $H(F_v)$ is the group of units   
of a central division algebra $D$ over $F_v$ of rank $n^2$. Let $\pi$ be
an irreducible, finite dimensional, complex smooth representation of $D^*$.
\ps

\begin{lemme} \label{pi} There exists an irreducible automorphic
representation $\Pi$ of $H(\AAA_F)$ such that: \ps
(a) If $x \neq v,\, w$ is finite place, then $\Pi_x$ is unramified, \ps
(b) $\Pi_v \simeq \pi\otimes \psi$ for some unramified character
$\psi: D^* \rightarrow \C^*$, \ps
\end{lemme}

\pf  We first choose, for each finite place $x$ of $F$, a particular compact
open subgroup $J_x$ of $H(F_x)$. If $x$ is different from $v$ and $w$,
$H(F_x)$ is quasisplit so that we can take for $J_x$ a maximal compact
subgroup which is {\it very special} in the sense of \cite[\S 3.6]{L} (when
$x$ does not ramify in $E$, the hyperspecial compact
subgroups of $H(F_x)$ are very special; for almost all $x$ we may and want to take $J_x=H(\OO_{F_x})$). For such a place $x$, an
irreducible admissible representation of $H(F_x)$
will be said to be {\it unramified} if it has a nonzero vector invariant by
$J_x$. If $x=v$, we take $J_x={{\mathcal O}_D}^*$ and we
fix an irreducible constituent $\tau_v$ of $\pi_{|J_v}$. If $x=w$, 
we take any compact open subgroup of $H(F_x)$ for $J_x$. Let $J:=\prod_x
J_x$, it is a
compact open subgroup of
$H(\AAA_{F,f})$. Let $\tau$ be the trivial extension to $J$ of the
representation $\tau_v$ of $J_v$, via the canonical projection $J
\rightarrow J_v$. \par  
        As $H_{\infty}:=H(F \otimes_{\Q} \R)$ is compact, the group  
$H(F)$ is discrete in $H(\AAA_{F,f})$. In particular, $\Gamma:=H(F) \cap J$
is a finite group. For any continuous (finite dimensional), complex
representation $W$ of $H_{\infty}$, we set $W(\tau):=W \otimes \tau^*$, viewed as
an $H_{\infty} \times J$-representation. By the Peter-Weyl theorem, we can
find an irreducible $W$ such that $W{|\Gamma}$ contains a copy of $\tau$, and so
a nonzero element $v \in W(\tau)^{\Gamma}$ (for the diagonal action of $\Gamma$). We choose
moreover an element $\varphi \in W(\tau)^*$ such that $\varphi(v)=1$. Let $h: H_{\infty} \times J \rightarrow
\C$ be the coefficient of $W(\tau)$ defined by $h(z):=\varphi(z^{-1}.v)$. By construction, $h$ is smooth, left
$\Gamma$-invariant, and generates copies of $W(\tau)^*$ under the right translations by $H_{\infty}\times J$,
since $W(\tau)$ is irreducible. It is nonzero since $h(1)=1$. It extends then
uniquely to a smooth map $f_h: H(\AAA_F) \rightarrow \C$ null outside the open subset
$H(F).(H_{\infty}\times J)$ and satisfying $f_h(\gamma z)=f_h(z), \, \, \forall \gamma \in H(F), \, z \in
H(\AAA_F)$, i.e. which is an automorphic form for $H$. Let $\Pi \subset
L^2(H(F)\backslash H(\AAA_F),\C)$ be an irreducible constituent of the
$H(\AAA_F)$-representation generated by $f_h$. \par By definition, $\Pi_x$ is
unramified if $x \neq v, w$ is a finite place, hence (a) holds. Moreover, $\Pi_v$ is an
irreducible representation of $H(F_v)$ whose restriction to $J_v$ contains $\tau_v$. 
As $\pi$ is supercuspidal and
as $(J_v,\tau_v)$ is a  $[H(F_v),\pi]_{H(F_v)}$-type by \cite[prop. 5.4]{BK}, 
$\pi$ and $\Pi_v$ can differ only by a twist by an unramified character of $H(F_v)$, proving (b).
We recall the argument for the convenience of the reader. Let $Z=F_v^*$ denotes
the center of $D^*$, $\widetilde{J_v}:=Z\mathcal{O}_D^* \subset D^*$, 
and $\widetilde{\tau}$ the natural extension of $\tau$ to $\widetilde{J_v}$ which is contained in
${\Pi_v}_{|\widetilde{J_v}}$. Note that any unramified character of $Z$ extends to a
unramified character of $D^*$, hence we may assume that
${\pi}_{|\widetilde{J_v}}$ contains $\widetilde{\tau}$, by replacing $\pi$
by some unramified twist if necessary. Note that $D^*/\widetilde{J_v} \simeq
\Z/n\Z$ is finite abelian and that its characters are unramified characters of $D^*$. We
conclude as $\Pi_v$ is
a constituent of $\Ind_{\widetilde{J_v}}^{D^*} \pi_{|\widetilde{J_v}} \simeq
\pi \otimes_{\C} \C[D^*/\widetilde{J_v}]$.\epf

\section{Proof of theorem I} \label{proof}

\subsection{Construction of $S$-unramified number fields} \label{HT} With the assumptions of \S \ref{constructionI}, let $S$
be the set of places of $E$ dividing  $lv$. We fix an
embedding $\varphi: E_S \rightarrow \overline{F_v}$ extending the
$F$-embedding $E \rightarrow
\overline{F_v}$
given by $u$. Let us denote by $$\G_{E,S}:=\Gal(E_S/E).$$
Attached to $\varphi$ is a group homomorphism $\G_{F_v} \rightarrow \G_{E,S}$ (see
the introduction for the notations). 
We keep the assumption of \S \ref{auto}, and
we choose a $\Pi$ given by lemma
\ref{pi}. We assume from now that {\it $\pi$ corresponds to a supercuspidal
representation of $\GL_n(F_v)$ by the Jacquet-Langlands correspondence},
that is $\dim_{\C}(\pi)>1$ if $n>1$. The local Langlands correspondence
(see \cite{HT}) associates to $\pi$ an continuous, irreducible, representation 
$$\psi_\pi: \W_{F_v} \longrightarrow \GL_n(\C).$$ 
Let $F_v^{ur}(\pi) \subset
\overline{F_v}$ be the finite
extension of $\Fur$ which is fixed
by $\Ker((\psi_\pi)_{|I_{F_v}})$, that is the extension of 
$F_v^{ur}$ cut out by $(\psi_{\pi})_{|I_{F_v}}$. Recall that $l$ is prime to $v$ and fix embeddings $\iota_{l}: \overline{\Q}
\rightarrow \overline{\Q_l}$ and $\iota_{\infty}: \overline{\Q} \rightarrow
\C$. 

\begin{lemme} \label{rep} (i) There exists a continuous representation $$R:
\G_{E,S} \rightarrow \mathrm{GL_n}(\overline{\Q_l}),$$ such that $R_{|W_{F_v}}$
corresponds
to $\Pi_v\otimes|\det|^{(n-1)/2}$ by the local Langlands correspondence, \ps

(ii) $\Fur.\varphi(E_S) \supset F_v^{ur}(\pi)$.\ps

\end{lemme}

\newcommand{\Res}{\mathrm{Res}}
\newcommand{\JL}{\mathrm{JL}}
\newcommand{\BC}{\mathrm{BC}}
\pf This lemma is a consequence of conditional automorphic base change
(\cite{C}, \cite{CL}, \cite{HL}), of Jacquet-Langlands correspondence, and of the main 
theorem of Harris \& Taylor \cite{HT}. Precisely, let us denote by $\BC$ the
quadratic base
change from $H$ to $ H':=\Res_{E/F}(H \times_F E)$. By theorem 3.1.3 of
\cite{HL}
(generalizing \cite[thm. A.5.2]{CL} and \cite{C}), there is a cuspidal
automorphic
representation $\Pi'$ of
$H'(\AAA_F)$ such that $\Pi'_x=\BC(\Pi_x)$ for each place $x \neq w$ of $F$.
It applies because $\Pi_{\infty}$ is automatically cohomological as $H(F\otimes_\Q\R)$ is compact, and because
$H$ is attached to a division algebra, as $H(F_v)$ is. Note that $\BC$ has
been defined for unramified representations at the places $x$
such that $H(F_x)$ is a ramified, quasisplit, unitary group, in \cite[\S
3.6, prop. 3.6.4]{L}\footnote{As
Labesse explained to us, the confusing hypothesis "$G$ is unramified, and
$K^G$ is hyperspecial" in prop. 3.6.4 {\it loc. cit.}
should be understood as "$G_0$ is unramified {\it over} $E_x$, and $K^G$ is
hyperspecial viewed as a subgroup of $G_0(E_x)$" (for us
$G_0=\GL_n(E_x)$), so as to be of any use, and the same proof applies
verbatim. If we didn't use these facts, we
would be obliged to add, in the set $S$ of the theorem, the places of $E$
ramified above $F$.
Note, however, that this would suffice to get corollary \ref{cor2}.}. In
particular, $\Pi'_x$ is
unramified if $x \notin S$, $\Pi'_v=\Pi_v\otimes \Pi_v^*$, and
$\Pi'_{\infty}$ has the base change infinitesimal character. We can then
apply\footnote{Note that $\JL(\Pi')$ is cuspidal by Moeglin-Waldspurger's
description of the discrete
spectrum of $\GL_n$, since
$\JL(\Pi')_v=\JL(\Pi'_v)$ is supercuspidal.} theorem C of \cite{HT} to the
image $\JL(\Pi')$ of $\Pi'$ by
the Jacquet-Langlands
correspondence (due to Vign\'eras, see \cite[thm. VI.1.1]{HT}), and consider
$$R:=R_l(\JL(\Pi'))$$ given by {\it loc. cit} and the embeddings $\iota_l,\,
\iota_{\infty}$, which proves (i). \ps
        We check the second assertion. Let us denote by $E_S(R)$ 
the subfield of $E_S$ fixed by $\Ker(R)$. As
$\Gal(\overline{F_v}/(F_v.\varphi(E_S(R))))=\Ker(R_{|\G_{F_v}})$, Galois theory shows that 
$$F_v.\varphi(E_S(R))=F_v(R_{|\G_{F_v}}).$$ \noindent By (i), $\psi_{\Pi_v}$ is defined over $\iota_{\infty}(\overline{\Q})$ and $R_{|I_{F_v}}\simeq
i_l.i_{\infty}^{-1}(\psi_{\Pi_v})_{|I_{F_v}}$ (it has finite image), so that
$\Fur(R_{|I_{F_v}})=\Fur(\psi_{\Pi_v})$. But by lemma \ref{pi}, $\pi$ and
$\Pi_v \otimes |\det|^{(n-1)/2}$ only differ by an unramified twist, hence
$\psi_{\Pi_v\otimes |\det|^{(n-1)/2}}$ and $\psi_{\pi}$
are equal when restricted to $I_{F_v}$. In particular, $\Fur(\pi):=\Fur(\psi_{\pi})=\Fur(\psi_{\Pi_v})$, and we conclude
as $\Fur.F_v(R_{|\G_{F_v}})=\Fur(R_{|I_{F_v}})$.\epf

\subsection{End of the proof} \label{trick} 

\begin{lemme} \label{sscorps} Let $\Q_p \subset M \subset L \subset  
\overline{M}$ be a tower of field extensions with $M/\Q_p$ finite.\ps
(i) Assume that $L/M$ is Galois and $M^{ur}.L=\overline{M}$, then $L=\overline{M}$. \ps
(ii) If $\sigma \in \G_M$ acts trivially by conjugation on
the tame inertia of $I_M$, then $\sigma \in I_M$.\ps
(iii) Assume that for all finite Galois extension $K/M$ such that $\Gal(K/M)$
admits an
injective, irreducible, complex linear representation, we have $K \subset
L$; then
$L=\overline{M}$. \ps
\end{lemme}

\pf (i) Let $H:=\Gal(\overline{M}/L)
\subset \G_M$. By assumption, $H$ is normal in $\Gamma$
and
$H\cap I_M=\{1\}$, hence $I_MH \simeq I_M \times H$ is a direct product 
and (i) (i.e. $H=\{1\}$) is an immediate consequence of (ii), which we now
prove. \ps
Let $I_M^{tr}$ be the tame
inertia quotient of $I_M$, recall that there is an isomorphism
$$I_M^{tr} \overset{\sim}{\longrightarrow}\prod_{l \neq p}
\Z_l(1),$$
which means that the action of $\Gamma$ by conjugation on $I_M^{tr}$, which
factors through the canonical map $\nu: \Gal(M^{ur}/M)
\longrightarrow \widehat{\Z}$, is the multiplication by
$q^{\nu(\cdot)} \in \prod_{l \neq p} {\Z_l}^*$, where $q$ is the cardinal of
the residue field of $K$.
Let $\gamma \in \Gamma$ such that $q^{\nu(\gamma)}=1$, we aim to prove
that $\nu(\gamma)=0$. 
We claim that for any integers $m \geq 2$ and $r\geq 1$, 
we can find infinitely many primes $l$ such that $r$ divides the
order of the image of $m$ in ${\mathbb{F}_l}^*$. Applying this claim to
$m=q$, we get that $r$ divides $\nu(\gamma)$ for all $r \geq 1$, i.e.
$\nu(\gamma)=0$.
Let us prove the claim now. Note that if $P \in \Z[X]$
is such that $P(0)\neq 0$ and $m\geq 2$ is an integer, then the sequence
$(P(m^{N!}))_{N \geq 1}$ admits infinitely many prime
divisors. Indeed, assuming the assertion it is false, we can find a
prime $p$ such that each power of $p$ divides $P(m^{N!})$ for some $N$. As $P(0)\neq 0$, $p$ is
prime to $m$. Thus if we let $N$ goes to
infinity, we get that $P(1)=0$ in $\Z_p$. This is a contradiction, as $P=X-1$ certainly
satisfies the assertion. The claim follows then from the special case
$P:=$ the $r^{th}$ cyclotomic polynomial.\ps
        We now show (iii). Let $K/M$ be any finite
Galois extension, $\rho_1,\dots,\rho_t$ the irreducible, complex, linear
representations of $\Gal(K/M)$, and $K_i \subset K$ the fixed field of
$\Ker(\rho_i)$.
By assumption, $K_i \subset L$. The existence of the (faithful)
regular representation of $\Gal(K/M)$ implies that 
$\bigcap_i \Ker(\rho_i)=\{1\}$, so that by Galois theory we have
$K=K_1\dots K_t$, and $K \subset L$. \epf

\medskip

Let us finish the proof of theorem \ref{T1}. Let $K/F_v$ be a finite Galois extension 
such that $\Gal(K/F_v)$ admits an irreducible injective representation $\rho:
\Gal(K/F_v) \rightarrow \mathrm{GL}_n(\C)$. The local Langlands
correspondence and the
Jacquet-Langlands correspondence associate to $\rho$ 
an irreducible smooth representation $\pi$ of $D^*/F_v$ such that    
$\psi_{\pi}$ factors through $\rho$. Lemma \ref{rep} (ii) shows that 
$\Fur.\varphi(E_S) \supset F_v^{ur}(\pi)$. But $F_v^{ur}(\pi)=\Fur.K$, as
$\rho$ is injective. We have
thus shown 
that $K \subset \Fur.\varphi(E_S)$. Applying lemma \ref{sscorps}
(iii) to $M=F_v$ and $L=\Fur.\varphi(E_S)$, we conclude that
$\Fur.\varphi(E_S)=\overline{F_v}$. Applying 
now lemma \ref{sscorps} (i) to $M=F_v$ and $L=F_v.\varphi(E_S)$, we get 
theorem \ref{T1}. $\square$ \ps

\section{Construction of some automorphic forms II}
\label{constructionII}
\newcommand{\g}{{\got g}}
\renewcommand{\k}{{\got k}}
\renewcommand{\c}{{\got c}}
\newcommand{\p}{{\got p}}
\newcommand{\n}{{\got n}}
\renewcommand{\b}{{\got b}}
\newcommand{\h}{{\got h}}
\newcommand{\cal}{\mathcal}

In this section, we generalize some results of \S \ref{constructionI} concerning the 
construction of automorphic representations with prescribed properties. Let $H$ be a connected reductive group over $\Q$, $J$ be a compact open subgroup of
$H(\AAA_{\Q,f})$, and let $\tau$ be a fixed, irreducible, smooth representation of $J$. We aim to construct discrete, 
algebraic, automorphic representations $\Pi$ of $H$, such that ${\Pi_f}_{|J}$ contains $\tau$. From the point of 
view of types theory (see \cite{BK}), it allows us to prescribe the inertial equivalence class of each $\Pi_v$ with $v$ finite.
For example, in the case of $\GL_n$ over a local field, the inertial equivalent class of a smooth irreducible representation determines the 
restriction to the inertia group of its Langlands parameter, hence types control the ramification in a strong sense. 
In general, it seems hard to prescribe more properties of $\Pi_v$, as we certainly cannot guess easily its Weil numbers if $v$ is unramified 
for instance. Note moreover that such a $\Pi$ may not exist, as the obstruction given by units shows in the example of the torus $H=\mathrm{Res}_{E/\Q}(\mathbb{G}_m/E)$
when $E \neq \Q$ is not quadratic imaginary. 
 
\subsection{} Let $H$ be as above and let $Z$ denotes its center. We assume that:\ps

(A) $H(\R)$ has a holomorphic discrete series (i.e. the derived Lie algebra of 
$H(\R)$ satisfies (\ref{hyphc}) of \S \ref{HC}), \ps

(B) $Z$ is discrete in $Z(\AAA_{\Q,f})$.\ps

\begin{prop} \label{general} Let $H$ be as above, $J$ be a compact open subgroup of
$H(\AAA_{\Q,f})$, and let $\tau$ be an irreducible smooth representation of $J$. 
There exists an irreducible, discrete, automorphic representation $\Pi$ of $H(\AAA_{\Q})$ such
that $\Pi_{\infty}$ is in the holomorphic discrete series, and ${\Pi_f}_{|J}$ contains $\tau$.
\end{prop}

{\it Remarks :} Note that by \cite[thm. 8.9]{B1}, property (B) depends only of the isogeny
class of $Z$. If $Z$ is anisotropic, it holds if and only if $Z(\R)$ is
compact by thm. 8.7. {\it loc.cit.} It happens for example when $E/F$ is CM and 
$Z$ is $(\mathrm{Res}_{E/\Q} \mathbb{G}_m)^{N_{E/F}=1}$. It is easy to see that (A) and (B) both hold when $H$ is $\GL_2/\Q$, $\GSp_{2n}/\Q$, 
the scalar restriction to $\Q$ of a unitary group attached to a CM field, $\mathrm{SO}(n,2)/\Q$, or when $H(\R)$ is compact. \ps

Before proving proposition \ref{general}, we mention its following
corollary. Let $E$ be a CM field, $F$ its maximal totally real subfield, $v$ a finite place of $F$ which splits in $E$, and write $v=uu'$.

\begin{cor}  \label{corgeneral} Let $\pi$ be a supercuspidal, irreducible, representation of
$\GL_n(E_u)$. There exists an automorphic cuspidal representation $\Pi$ of
$\GL_n(\AAA_E)$ satisfying $\Pi^c \simeq \Pi^*$, and such that: \ps

(a) $\Pi$ is unramified outside $u, u'$, and $\Pi_u$ is an unramified twist
of $\pi$, \ps

(b) for any $x$ archimedian, $\Pi_x$ has a regular algebraic infinitesimal character. 
\end{cor}

\pf By the same reasoning as in the proof of lemma \ref{groupe}, Hasse's principle shows that there exists a unitary group 
$\U(n)/F$ attached to $E/F$ which is quasisplit at all finite places $\neq v$, and the group of invertible elements of a central division 
algebra over $F_v \overset{\sim}{\rightarrow} E_u$ at $v$. In fact, as the local invariant of $\UU(m+r,m-r)$ is 
$r \bmod 2$ by \cite[lemme 2.2]{C}, we could even assume that $\U(n)$ is compact at all archimedian places, except maybe one where it 
is $\UU(n-1,1)$. Then $H:=\mathrm{Res}_{E/\Q} \U(n)$ satisfies the hypothesis of proposition \ref{general}, by the remarks above.
Let $\pi':=\JL^{-1}(\pi)$ be the finite dimensional representation of
$\U(n)(F_v)$ corresponding to $\pi$ by the Jacquet-Langlands correspondence. 
By choosing a suitable compact open subgroup $J$ of $H(\AAA_{\Q,f})$ as in the proof of lemma \ref{pi}, proposition \ref{general} produces an 
automorphic discrete representation $\Pi'$ of $U(n)/F$ which is unramified outside $v$, such that $\Pi'_v$ is
a twist of $\pi'$ by an unramified character, and whose archimedian components are discrete series (hence cohomological by \cite[thm. 5.3. (b)]{BW}). As we are in the 
hypothesis of the ``conditional'' quadratic base change, we conclude as in the proof of lemma \ref{rep} (i) that $\Pi:=\JL(\BC(\Pi'))$ satisfies all the 
hypothesis of the corollary. \epf

\ps
Note that proposition \ref{general} is a generalization of lemma \ref{pi}, and the same proof as {\it loc. cit.} would allow us
to conclude verbatim when $H(\R)$ is compact\footnote{For the sake of
exposition, and as much of the difficulty disappears in the compact case, we found it convenient to separate the proof. Moreover, lemma \ref{pi}
is not only sufficient for the proof of theorem \ref{T1}, but also fix some ideas about the general
case.}. The general case relies on more delicate analytic facts. We are very grateful to Laurent Clozel 
for explaining us an argument using the Selberg trace formula, holding even
under the assumption that $H(\R)$ has a discrete series. 
In what follows, we will give an argument using Poincar\'e series
going back to Poincar\'e, H. Cartan and Godement (see \cite[ch. 6]{B2}, \cite{Ca1} \S 1, 
\cite{Ca2} \S 6, \S 10, and \S 10 bis, for the detailed case
$H_{\infty}=\mathrm{GSp}_{2g}(\R)$, and most of the ideas of the general argument), so that proposition \ref{general} should certainly not be considered as original. Roughly, the difficulty is twofold. First, we want to use coefficients of
square integrable representations of $H(\R)$ to ensure the convergence of some Poincar\'e
series. Then, we must have a sufficiently good control on these coefficients to produce some non identically zero automorphic forms. 
Thus, as a preliminary, will recall in \S\ref{HC} some facts about Harish-Chandra's realization of the holomorphic discrete series (due to Harish-Chandra). \ps

\newcommand{\ad}{\mathrm{ad}}

\subsection{Holomorphic discrete series} \label{HC} In this section, we will
follow closely Knapp's book \cite[chap. VI]{K}. Let $G$ be a real, connected, reductive
Lie group, and $K$ a maximal compact subgroup. Assume that a Cartan decomposition $\g=\k\oplus \p$ satisfies 
\begin{equation} \label{hyphc} Z_{\g}(\c)=\k, \end{equation} \noindent where $\c$ is the center of $\k$.
For example, it is easily seen 
to be the case for any $\UU(p,q)$ and for symplectic groups (see {\it loc. cit. } \S
2). Let $\h \subset \k$ be a Cartan subalgebra, associated to a Cartan subgroup
$T \subset K$, and let $\g^{\C}=\h^{\C}\oplus_{\alpha \in \Delta} \g_{\alpha}$ be its associated root space
decomposition. \ps 
        Recall that a root $\alpha \in \Delta$ is said to be {\it noncompact} (resp. {\it compact}), if $\g_{\alpha} \subset \p^{\C}$ (resp. in $\k^{\C}$). This
notion gives us a partition $\Delta=\Delta_n \coprod \Delta_K$. We fix a
{\it good} ordering on $\h_{\R}^*$, i.e. with the property that any positive noncompact root is
bigger than any compact root. We set $\p^{\pm}:=\oplus_{\alpha \in
\Delta_n^{\pm}} \g_{\alpha}$, these are abelian Lie subalgebras of $\p^{\C}$
which are stable by $\ad(\k)$. We fix a complex matrix group $G_{\C}$ whose
Lie algebra is $\g^{\C}$, and we let $K_{\C}$ and $P^{\pm}$ be the complex
analytic subgroups of $G_{\C}$ with Lie algebras $\k^{\C}$ and $\p^{\pm}$, 
as in \S 3 {\it loc.cit}. By Harish-Chandra's
decomposition, the natural product map $P^+ \times K_{\C} \times P^-
\rightarrow G_{\C}$ is a complex open immersion, and $G K_{\C}P^- \subset
G_{\C}$ is an open subset, and so inherits the complex structure of
$G_{\C}$. Moreover, there exists a bounded open domain $\Omega \subset
P^+\simeq \C^p$ such that $\Omega K_{\C} P^-=G K_{\C}P^-$. As $G \cap K_{\C} P^-=K$,
the $G$-orbit of $1 \in G_{\C}/K_{\C}P^-$ is $G/K\simeq \Omega$. \ps

Let us fix an analytically integral $\lambda \in
\h_{\R}^*$, which is dominant restricted to $\Delta_K^{+}$, and let $W_{\lambda}$
be the irreducible representation of $K_{\C}$ of highest weight $\lambda$.
This $W_{\lambda}$ gives rise to a complex algebraic vector bundle on the
grassmannian $G_{\C}/K_{\C}P^-$, and we want to consider holomorphic, square
integrable, sections of this bundle on $\Omega$. Concretely, we fix
a definite, $K$-invariant, hermitian product on $W_{\lambda}$, and we consider the vector space
$V_\lambda$ of $W_\lambda-$valued functions $f$ on $G$
satisfying:\ps
\begin{itemize}
\item[] (i) $\forall g \in G$, $\forall k \in K$, $f(gk)=k^{-1}.f(g)$,\ps
\item[] (ii) $f$ is holomorphic (see below),\ps
\item[] (iii) $\int_{G} |f(g)|^2dg < \infty$.\ps
\end{itemize}

For a map $f: G \rightarrow W_{\lambda}$ satisfying (i), (ii) means that
the canonical $P^-$-invariant extension of $f$ to $GK_{\C}P^-$ is holomorphic. 
Note that $V_{\lambda}$ is in a natural way a module over the ring of bounded holomorphic functions on
$\Omega$, and a unitary representation $L_{\lambda}$ of $G$ by left translation. 
Let $\delta:=\frac{1}{2}\sum_{\alpha \in \Delta^{+}} \alpha$. We can now state Harish-Chandra's theorem 
(\cite[thm. 6.6]{K}\footnote{The identification of the space we call
$V_{\lambda}$ with the one {\it loc.cit.} comes from Borel-Weil-Bott
realization of $W_{\lambda}$ \cite[thm. 5.29]{K}.}):\ps
%if $f(g)$ is viewed as a function on $K_{\C}$ such that
%$f(kt)=\chi_{\lambda}(t)^{-1}f(k)$, then $f\mapsto f'$ is an intertwining
%if we set $f'(g)=(f(g))(1)$.
{\bf Theorem:} (Harish-Chandra) Assume that $\langle\lambda+\rho,\alpha \rangle<0$ for all $\alpha
\in \Delta_n^+$. Then $V_\lambda$ is a non zero Hilbert space, and $L_{\lambda}$ is
an irreducible, square-integrable, $G$-representation. \ps

\smallskip

        We consider now a special case related to the
standard jacobian automorphy factor on $\Omega$. Let $\rho_n:=\sum_{\alpha
\in \Delta_n^+} \alpha  \in \h_{\R}^*$. We claim that $\rho_n$ is 
analytically integral, dominant, and satisfies:
\begin{equation} \label{pos} \langle \rho_n , \alpha \rangle > 0,  \, \, \, \, \forall
\alpha \in \Delta_n^+ \end{equation}

Indeed, if $p:=\dim_{\C}\p^+$, then $\rho_n$ is the weight of the $1$-dimensional
representation $\Lambda^p \p^+$ of $\k$. In particular, $\forall m \in
\Z$, $m \rho_n$ is
analytically integral, and dominant with respect to $\Delta_K^+$. Moreover, by definition of the ordering, 
$\rho_n$ is also the highest weight of the natural representation of $\g$ on
$\Lambda^p \g^{\C} \supset \Lambda^p \p^+$, and so is dominant. The last part of the claim
follows, as the parabolic subalgebra $\k^{\C} \oplus
\p^+$ is its own normalizer in $\g^{\C}$. \ps
        As a consequence, we can fix an
integer $r>0$ such that $\delta:=-r \rho_n$ satisfies the hypothesis of
Harish-Chandra's theorem. Let us fix a an element 
$j \in V_{\delta}$ which is an eigenvector for the left translations by $K$,
and satisfies $|j(1)|=1$. Such an element exists: we can take the one denoted
by $\psi_{\lambda}$ in \cite[lemma 6.7]{K}. By Harish-Chandra's
convolution theorem (\cite[cor. 8.41]{K}, \cite[cor. 2.22]{B2}), $K$-finite elements of
$V_{\lambda}$ are bounded on $G$, and so is $j$. In particular, if $f \in
V_{\lambda}$ and if $m\geq 0$ is an integer, then $(f\cdot j^m)(g):=f(g)
\otimes j^m(g)$ defines an element in $V_{\lambda+m\delta}$. 
For convenience, we shall often identify the $1$-dimensional vector space $W_{\delta}$ with $\C$, taking $j(1)$ as
a norm $1$ basis element.\ps

{\it Example:} \, \, Let $G:=\UU(n-1,1)$, let $K$ be the diagonal $\UU(n-1) \times \UU(1)$, and let $T \subset K$
be the diagonal torus. We can choose a good ordering on
the roots such that $\sum_{\alpha \in \Delta^+} \g_{\alpha}$ is the standard
upper Borel subalgebra in $\g^{\C}=\mathrm{gl}_n(\C)$. In this case,
$K_{\C}P^-$ identifies with the standard lower parabolic of $G_{\C}=\GL_n(\C)$ of
type $(n-1,1)$. The bounded domain $\Omega$ identifies with
the open unit hermitian ball of $\C^{n-1}=P^+ \subset
G_{\C}/K_{\C}P^-=\mathbb{P}^{n-1}(\C)$. We define $e_i \in \h_{\R}^*$ by $e_i(x_1,\dots,x_n):=x_i$. 
With these choices, $\Delta_n^+=\{e_i-e_n, 1\leq i <n\}$
and $\Delta_K^+=\{e_i-e_j, \, \, 1\leq i<j<n\}$. 
Then $\lambda=\sum_{i=1}^n m_i e_i$ satisfies the hypothesis of Harish-Chandra's theorem
if, and only if, $m_1 \geq m_2 \geq \cdots \geq m_{n-1}$, $m_n > m_1+(n-1)$,
and $\forall i, m_i \in \Z$.
We have $\rho_n=(e_1+\dots+e_{n-1})-(n-1)e_n$, and we can take $\delta=-2\rho_n$.
\ps

\subsection{Convergence of Poincar\'e series} We now begin the proof of
proposition \ref{general}. Let $Z_c$ denotes the maximal compact subgroup of $Z(\R)$ (viewed as a real Lie group) 
and write $Z(\R)=Z_n \times Z_c$ for some closed subgroup $Z_n$ isomorphic to some $\R^m$.
Let $G$ denotes the closed connected subgroup of the Lie group $H(\R)$
whose Lie algebra is $\got{g}=[\got{h},\got{h}]\oplus \got{z_c}$. If $H(\R)^+$ denotes the neutral component of the real Lie group $H(\R)$, then $$H(\R)^+=Z_n\times G.$$
Let $\Gamma \subset H(\R)^+$ be the discrete subgroup $H(\Q) \cap (H(\R)^+\times J)$, let us show that $\Gamma \subset G$. 
Let $C$ be the $\Q$-torus $H/H_{der}$ and $\det: H \longrightarrow C$ be the canonical $\Q$-morphism. Note that $\det$ induces a $\Q$-isogeny $Z \longrightarrow C$, hence $C$ has finite arithmetic groups by assumption (B) and \cite[thm. 8.9]{B1}. In particular, $\det(\Gamma) \subset C(\Q)$ is a finite subgroup. But the induced map $\det: H(\R)^+ \longrightarrow C(\R)$ is injective restricted to $Z_n$, and $\det(G)$ is compact, hence $\Gamma \subset G$.  \par
The Lie algebra $\got{g}$ satisfies (\ref{hyphc}) by assumption (A), so that we can apply to $G$ the constructions (and notations) of \S \ref{HC}. Let $\lambda \in \h_{\R}^*$ be as in Harish-Chandra's theorem, we set $V_{\lambda}(\tau):=V_{\lambda}\otimes
\tau^*$, viewed as a representation of $G \times J$. We will consider
coefficients of $V_{\lambda}(\tau)$ of the following kind. Let 
$W_{\lambda}(\tau):=W_{\lambda}\otimes_{\C} \tau^*$ viewed as before as a
representation of $K \times J$. If $\varphi \in W_{\lambda}(\tau)^*$, we can
see it as a continuous linear form $\tilde\varphi$ on $V_{\lambda}(\tau)$ by $f \mapsto
\varphi(f(1))$. We define $h:=h_{f,\varphi}: G \times J \rightarrow \C$
by $h(z):=\tilde\varphi(z^{-1}.f)=\varphi(( 1 \times z_J^{-1}).f(z_G))$. The Poincar\'e 
series $P_h:G\times J \rightarrow \C$ is defined by

$$P_h(z) :=\sum_{\gamma \in \Gamma} h(\gamma z)$$

\noindent We equip $W(\tau)$ with any hermitian norm $|.|$ such that $K \times J$ acts by
unitary transformations. As $f \in L^2(G)$, it comes that $P_h \in L^2_{loc}(G
\times J)$. Indeed, if $U \subset G \times J$ is a compact set,
then $\Gamma \cap (UU^{-1})$ is finite, hence 
\begin{equation} \label{l2} \sum_{\gamma \in \Gamma} \int_{U} |h(\gamma g)|^2dg \leq
|\Gamma\cap UU^{-1}|\int_{G \times J}|h(g)|^2dg \leq C \int_G |f(g)|^2 dg <
\infty,
\end{equation}
where $C:={\rm vol}(I)\cdot |\Gamma\cap UU^{-1}| \cdot||\varphi|| >0$. As all the $g \mapsto f(\gamma.g)$ are in $V_{\lambda}$, and in particular
holomorphic, $P_h$ converges in fact {\it uniformly} on any compact subset of $G
\times J$, to a
holomorphic function on $G \times J$ of right $K \times J$-type $W(\tau)^*$. In
particular, we proved the

\begin{lemme} \label{conv} Let $\lambda \in \h_{\R}^*$ be as in
Harish-Chandra's theorem, and let $h$ be a coefficient of $V_{\lambda}(\tau)$ as
above. The Poincar\'e s\'eries $P_h$ is normally convergent
on any compact subset of $G \times J$. 
\end{lemme}

\subsection{Construction of non vanishing Poincar\'e series} It remains to
find a $\lambda$ and an $h$ such that $P_h \neq 0$. As in the proof of 
lemma \ref{pi}, we can find a finite
dimensional, irreducible, complex representation $W=W_{\lambda}$ of $K$ such that
$W_{\lambda}(\tau)^{\Gamma_0}\neq 0$. By twisting $W_{\lambda}$ by $W_{\delta}^m$ for an
integer $m\geq 0$ big enough and divisible by $|\Gamma_0|$, (\ref{pos})
shows that we may assume
that $\lambda \in \h_{\R}^*$ satisfies
the hypothesis of Harish-Chandra's theorem. We fix as before an element $v
\in W_{\lambda}(\tau)^{\Gamma_0}$ of norm $1$ and we choose a $\varphi \in W(\tau)^*$
of norm $1$ such that
$\varphi(v)=1$. By lemma \ref{conv} applied to $\lambda=\delta$ and
to the coefficients associated to the function $j$ itself, we get that
$$\Gamma_1:=\{\gamma \in \Gamma, |j(\gamma)|\geq 1\}$$ is a finite
subset of $\Gamma$, and $\Gamma_0 \subset \Gamma_1$. We claim that we can choose an
$f \in V_{\lambda}(\tau)$ such that
$$f(1)=v, \, \, \, \, {\rm and} \,\, \, \, \forall \gamma \in \Gamma_1\backslash \Gamma_0, \, \, \,
f(\gamma)=0.$$
As $V_{\lambda}(\tau)$ is nonzero and stable by the action of $G \times
J$, we can find a $f_1 \in V_{\lambda}(\tau)$, such that $f_1(1)\neq 0$.
But by irreducibility of $W_{\lambda}(\tau)$ as $K \times J$-representation, we can find in
$\C[K \times J].f_1$ an element sending $1$ to $v$. We can thus assume that
$f_1(1)=v$. Then, we can multiply $f_1$ by a bounded holomorphic function $f_2$ on
$\Omega$ such that
$f_2(1)=1$ and $f_2(\gamma)=0$ if $\gamma \in \Gamma_1 \backslash \Gamma_0$,
which certainly exists (e.g. restrictions to the
bounded symmetric domain of polynomials on $P^+=\C^p$). The function $f=f_1f_2$ does the trick. This
proves the claim. \ps

For any $m\geq 1$ divisible by $|\Gamma_0|$, let us consider the element
$v_m:=v\otimes 1^{\otimes m}\in W_{\lambda+m\delta}(\tau)^{\Gamma_0}$ 
and $\varphi_m:=\varphi \otimes {\rm id}^{\otimes m} \in
W_{\lambda+m\delta}(\tau)^*$. We still have $\varphi_m(v_m)=1$. We
define $h_m$ to be the coefficient associated to $f\cdot j^m$ and
$\varphi_m$. By all the previous choices, $|h_m(\gamma)| \leq
|h_0(\gamma)|, \, \, \forall \gamma \in \Gamma$. As $\sum_{\gamma \in \Gamma} |h_0(\gamma)| < \infty$ by
lemma \ref{conv}, we get:

\begin{equation}\label{lim} P_{h_m}(1) \underset{|\Gamma_0|\,  | m, \, \, m \rightarrow \infty}\longrightarrow
\sum_{\gamma \in \Gamma_0} \varphi(\gamma.v)=|\Gamma_0| \neq 0 
\end{equation}
\noindent By (\ref{lim}), we 
may choose an integer $m$ big enough such that $P_{h_m}$ is not identically
zero. We denote again by $P_{h_m}$ its canonical $Z_n$-invariant
extension to $H(\R)^+ \times J=G\cdot Z_n \times J$. \ps
Let $f: H(\AAA_{\Q}) \rightarrow \C$ be the unique map which is
$H(\Q)$-invariant on the left, zero outside $H(\Q).(H(\R)^+\times J)$, and
which coincides with $P_{h_m}$ on the open subset $H(\R)^+\times J$ (note
that $f$ is well defined). Then 
$f$ is smooth, non zero, and belongs to $L^2(H(\Q)\backslash H(\AAA_{\Q})/Z_n,\C)$ by the estimate (\ref{l2})
(take $U$ a measurable fundamental domain for $\Gamma$ acting on $G \times
J$, they have finite volume). The closure of the $H(\AAA_{\Q})$-subrepresentation generated by $f$ is then a finite sum 
of topologically irreducible representations, any irreducible constituent
$\Pi$ of which satisfies by construction the hypothesis of proposition 
\ref{general}. \epf 

\medskip

{\it Remark :} \label{loin} As shown by the proof above and by lemma \ref{walls}
below, we can even assume that the parameter $\lambda$ of $\Pi_{\infty}$ 
is as far from the walls of $\h_{\R}^*$ as we want. 

\begin{lemme} \label{walls} Let $\Gamma$ be a finite subgroup of a compact connected Lie
group $K$, and let $\tau$ be an irreducible complex representation of $\Gamma$. For any real $C>0$, we can find
an irreducible representation $V$ of $K$ such that $V_{|\Gamma}$
contains $\tau$, and such that the highest weight of $V$ has distance $\geq C$ from the walls.
\end{lemme}
\pf Using character formulas of Weyl and Kostant, we could even prove that for
a dominant weight $\lambda$ far enough from the walls, $(V_{\lambda})_{|\Gamma}$ contains $\tau$ if and only if
they have the same central character when restricted to $\Gamma$. We thank Y. Benoist for explaining us the following quick proof. We fix $\h
\subset \k$ a Cartan subalgebra, choose a Weyl chamber, and we denote by $V_{\lambda}$ the complex irreducible representation of $K$ of
highest weight $\lambda \in \h_{\R}^*$. By the Peter-Weyl theorem, we can find
dominant weights $\lambda_1, \dots, \lambda_r$ such that $W:=\oplus_{i=1}
V_{\lambda_i}$ contains, when restricted to $\Gamma$, each irreducible finite dimensional complex representation
of $\Gamma$. Let $\lambda$ be any dominant weight, then the highest weights
of $V_{\lambda} \otimes_{\C} W$ lie in the ball of $\h_{\R}^*$ centered in
$\lambda$ and of radius $\sup_{i=1}^r |\lambda_i|$ (see \cite[IV. \S
11.13]{K}). As $\tau \otimes_{\C} W^*_{|\Gamma}$ contains all
irreducible representations of $\Gamma$ by construction, this concludes the
proof. \epf

\section{The case $E=\Q$}\label{Q} 

The problematic in this section is the following: is it possible to reduce
property $P_{E,S,u}$ to some standard conjectures in the arithmetic theory
of automorphic forms? To fix the ideas, we restrict to the case $E=\Q$. 
As proposition \ref{general} does not apply to $\GL_n/\Q$ for $n>2$, we 
must use other classical groups, which introduces some obstructions on the $L$-parameters we can reach with them.
It turns out that admitting the (certainly believed!) hypothesis $1$ and $2$ below, we can prove that 
$P_{\Q,\{\infty,l,p\},p}$ holds for any prime numbers 
$l \neq p$, which is much stronger than corollary \ref{cor2}. We will end the section 
by discussing $P_{\Q,\{\infty,p\},p}$ for a prime number $p$, which is the most interesting case in some sense, as its truth 
would imply all the results of this paper.
\newcommand{\GSpin}{\mathrm{GSpin}}
\newcommand{\tr}{\mathrm{tr}}
\newcommand{\SO}{\mathrm{SO}}
\subsection{Case $S=\{\infty,l,p\},\, \,  p\neq l$ (proof of theorem \ref{T2})}\label{conj} In this section, we give one way (certainly among many others!) to deduce $P_{\Q,\{\infty,l,p\},p}$ from some expected properties of automorphic forms for the symplectic groups $\GSp_{2n}/\Q$ for all $n\geq 2$. The proof is very similar to that of 
theorem \ref{T1}, so that we will be rather sketchy on some points. Let $\Q_p^{ab} \subset \overline{\Q_p}$ be the maximal abelian extension of $\Q_p$, note that $\Q_p^{ab} \subset 
\Q_p^{ur}.\Q_{\{\infty,p\}}$. \ps

\pf (of theorem \ref{T2} of the introduction)\, \,  Let $K/\Q_p$ be a finite Galois extension whose Galois group admits an
injective irreducible representation $\rho: \Gal(K/\Q_p) \rightarrow
\GL_n(\C)$, we want to prove that $K \subset
\Q_p^{ur}.\Q_{\{\infty,p\}}$. The representation $\rho$ extends to an
irreducible $L$-parameter $$\psi: \W_{\Q_p} \rightarrow \GL_n(\C).$$
Note that $\Q_p^{ab} \subset
\Q_p^{ur}.\Q_{\{\infty,p\}}$ by local class field theory, hence we can 
replace $\psi$ by any of its twists $\psi'$ by a continuous character
$\W_{\Q_p} \rightarrow \C^*$. We claim that for some well chosen character,
the dual of $\psi'$ is not isomorphic to any unramified twist of
$\psi'$. Indeed, as ${\psi}_{|I_{\Q_p}}$ has a finite image, we can find a $g \in
I_{\Q_p}\cap\Ker(\psi)$ acting on $\Q_p(\mu_{p^\infty})$ as an element of
infinite order. We can then find a finite order character $\chi$ of
$\Gal(\Q_p(\mu_{p^{\infty}})/\Q_p)$ such that $\chi^2(g)\neq 1$. Hence
$\psi':=\psi \otimes \chi$ does the trick as
$\tr(\psi'(g))=n\chi(g) \neq n\chi(g)^{-1}=\tr(\psi'^*(g))$. We can
therefore assume that $\psi^*$ is not isomorphic to any unramified twist of
$\psi$. 

\ps
Recall now that $M:=\GL_n \times \GL_1$ is a Levi factor of the maximal parabolic subgroup of
$\GSp_{2n}$ stabilizing a maximal isotropic subspace (the so-called {\it
Siegel parabolic}). 
As the root datum of $\GL_n \times \GL_1$ is selfdual, $\GL_n(\C)\times
\GL_1(\C)$ is also a Levi 
subgroup of $\widehat{\GSp_{2n}}=\GSpin_{2n+1}(\C)$. For any unramified
characters $\alpha,\beta: \W_{\Q_p} \rightarrow \C^*$, we can thus consider the
$L$-parameter $$\psi_{\alpha,\beta}:=\psi\otimes \alpha \times
\beta:\, \, \, \, \,  \W_{\Q_p}
\rightarrow \GSpin_{2n+1}(\C),$$ 
deduced by functoriality. As $\psi^*$ is not isomorphic to any  
unramified twist of $\psi$, the centralizer of any $\psi_{\alpha,\beta}$ in $\GSpin_{2n+1}(\C)$
is reduced to the center $\C^*\times \C^*$ of $\GL_n(\C) \times \GL_1(\C)$ (it obviously contains it). Indeed, if an element $\gamma \in \GSpin_{2n+1}(\C)$ satisfies $\gamma
\psi_{\alpha,\beta}(g)\gamma^{-1}=\psi_{\alpha,\beta}(g)$ for all $g \in
\W_{\Q_p}$, its image $\overline{\gamma}$ in $\mathrm{PGSpin}_{2n+1}(\C)=\SO_{2n+1}(\C)$
self-intertwines $$(\psi \otimes \alpha) \oplus (\psi\otimes
\alpha)^*
\oplus 1,$$ \noindent and we conclude because this semi-simple representation is
multiplicity free by assumption on $\psi$. As a consequence, we expect that the hypothetical $L$-packet of $\GSp_{2n}(\Q_p)$ associated to
any $\psi_{\alpha,\beta}$ have only one element. Namely, if
$\pi_{\alpha,\beta}$ denotes the supercuspidal representation of
$M=\GL_n(\Q_p) \times \GL_1(\Q_p)$ attached to $\psi \otimes \alpha \times \beta$ by
Harris \& Taylor, then the L-packet of $\psi_{\alpha,\beta}$ should consist of the full
\footnote{This induced representation is known to be irreducible, see \cite{Roc}.} 
normalized parabolic induction from $M$ to $\GSp_{2n}(\Q_p)$ of $\pi_{\alpha,\beta}$. We come now to our first hypothesis: \ps

{\sc Hypothesis 1:} For all supercuspidal representations $\pi$ of
$\GL_n(\Q_p)$ with no selfdual unramified twist, the inertial equivalence class 
$[M,\pi \times 1]_{\GSp_{2n}(\Q_p)}$ admits a type in the sense of \cite{BK}.\ps

By the above hypothesis applied to $\pi \times 1$, and by proposition \ref{general}, we can find
a discrete, irreducible, automorphic representation $\Pi$ of
$\GSp_{2n}(\AAA_{\Q})$ such that $\Pi_l$ is unramified if $l\notin\{\infty,p\}$, $\Pi_p \in [M,\pi \times 1]_{\GSp_{2n}(\Q_p)}$, and $\Pi_{\infty}$ belongs to the holomorphic discrete series. Moreover, by
remark \ref{loin}, we may assume that the parameter $\lambda$ of $\Pi_{\infty}$ is far from the
walls, hence $\Pi$ should conjecturally belong to a tempered $A$-packet.
Moreover, by the analysis preceding Hypothesis 1, the $L$-parameter $\Psi_p$ of $\Pi_p$ should
have the form $\psi_{\alpha,\beta}$ for some $\alpha, \beta$. 
By replacing $\Pi$ by an unramified twist if necessary, we may assume that $\Pi_{\infty}$ is algebraic. Let us fix a
complex injective representation $r: \GSpin_{2n+1}(\C) \rightarrow
\GL_{n_r}(\C)$, and embeddings $\iota_l$, $\iota_{\infty}$ as in \S \ref{HT}.
The reasoning above give credits to the following:\ps

{\sc Hypothesis 2:} For some $r$ as above, there exists a continuous representation 
$$R: \Gal(\Q_{\{\infty,p,l\}}/\Q) \rightarrow \GL_{n_r}(\overline{\Q_l}),$$
such that the representation $\iota_{\infty}.\iota_l^{-1}.R_{|I_{\Q_p}}
\simeq r \rond {\Psi_p}_{|I_{\Q_p}}$.\ps

Of course, following Langlands, the hypothetical Galois representation $R$ associated to
the automorphic representation $\Pi$, $\iota_{\infty}.\iota_l^{-1}$ and $r$,
should satisfy this hypothesis. Assuming the above hypothesis, we obtain by the same reasoning as
in the proof of lemma \ref{rep} (ii) that $\Q_p^{ur}.\Q_{\{\infty,p,l\}} \supset K$, which concludes the proof. \epf
\ps
{\it Remarks:} (i) Bushnell and Kutzko show in \cite{BK} \S 8 that Hypothesis $1$ would follow from 
the existence of $\GSp_{2n}(\Q_p)$-covers for the types they have constructed for the supercuspidal representations of $\GL_n(\Q_p)$. 
Some of theses covers have already been constructed by Blondel in \cite{B}, but
unfortunately only when $\pi$ has a selfdual unramified twist. \ps
(ii) Hypothesis $2$ is known for $n=2$ (Carayol, Deligne, Langlands). \ps

\subsection{Case $l=p$}\label{l=p} We would like to apply the same reasoning as in \S \ref{conj} to prove 
property $P_{\Q,\{\infty,p\},p}$, by choosing a $K$ as above and using $\Pi$, and using this time 
the hypothetical $p$-adic Galois representation $R$ attached to $\Pi$. 
The proof of \S \ref{conj} applies verbatim until the statement of hypothesis $2$, whose 
correct form becomes ($\Pi$ and $r$ are chosen as above): \ps

{\sc Hypothesis 2':} There exists a continuous representation 
$$R: \Gal(\Q_{\{\infty,p\}}/\Q) \rightarrow \GL_{n_r}(\overline{\Q_p}),$$
which is potentially semistable at $p$, and such that
$\iota_{\infty}.\iota_p^{-1}.(R_{|\G_{\Q_p}})^W_{|I_{\Q_p}}
\simeq r \rond {\Psi_p}_{|I_{\Q_p}}$.\ps

\noindent In the above statement, $\rho^W$ is the representation of
$\W_{\Q_p}$
attached by Fontaine to a $p$-adic $pst$ ({\it i.e.} potentially semistable) representation $\rho$ of
$\G_{\Q_p}$ (in fact, in the special situation above, the monodromy operator $N$ is automatically $0$). 
Let $\rho:=R_{|\G_{\Q_p}}$. The last step to conclude would be to establish a link 
between $F_1:=\Q_p^{ur}(\rho_{|I_{\Q_p}})$ (which is 
also $\Q_p^{ur}.\varphi(\Q_{\{\infty,p\}}(R))$) and $F_2:=\Q_p^{ur}(\rho^W_{|I_{\Q_p}})$ (which is
$\Q_p^{ur}.K$ in the notation of the proof of \ref{conj}). 
Note that $F_1/\Q_p^{ur}$ (resp. $F_2$) is in general infinite (resp. always finite), it is the smallest algebraic extension (resp. finite Galois extension) 
of $\Q_p^{ur}$ over which $\rho$ becomes trivial (resp. semistable, even crystalline here). It turns out that, it is not always the case that 
$F_2 \subset F_1$, which enables us to conclude as before. We discuss this point in the following subsection. \ps

{\it Remark:} The "modular form case" of hypothesis $2'$
(i.e. $n=2$) is known by the work of Carayol, Deligne, Langlands and Takeshi Saito (see the 
main theorem in \cite{S}). In general, hypothesis $2'$
would follow from hypothesis $2$, the construction of the pure motive attached to $\Pi$ and
 the conjectured {\it independence of $l$} of the
semisimplified Weil-Deligne representation
associated to the $l$-adic cohomology of a given proper smooth scheme $X$
over $\Q_p$ (see \cite{F} \S 2.4.3). \ps

\subsection{Kernels of $pst$ representations} Let $l\neq p$ be a prime number, we recall first the $l$-adic situation. 
Let $M/\Q_p$ be a local field. We fix a $\varphi  \in \W_M$ lifting a geometric 
Frobenius element, and a non-zero continuous homomorphism $t_l: I_M \rightarrow \Q_l$.  Let $\rho: \G_M
\rightarrow \GL_d(\overline{\Q_l})$ be a continuous representation. By
Grothendieck's $l$-adic monodromy theorem, there exists
a unique nilpotent matrix $N \in \mathrm{M}_d(\overline{\Q_l})$ such that $\rho$ and the representation of $I_M$
defined by $\gamma \rightarrow \exp(t_l(\gamma)N)$ coincide on some open subgroup of $I_M$. The Weil-Deligne representation attached 
to $\rho$ is then the isomorphism class of the pair $(\rho^W,N)$ where $\rho^W: \W_M \rightarrow \GL_d(\overline{\Q_l})$ is defined by 
(see \cite{Ta} \S 4.2, \cite{F})
$$ \rho^W(\varphi^n \gamma):=\rho(\varphi^n \gamma) \exp(-t_l(\gamma)N), \,\,\,\,\,\,
\gamma\in I_M,\,\,n\in \Z.$$
Note that $\rho^W_{|I_M}$ has finite image, hence consists of semisimple elements. In particular, $\Ker(\rho_{|I_M}) \subset 
\Ker(\rho^W_{|I_M})$, that is 
\begin{equation}\label{inclusion} M^{ur}(\rho^W_{|I_M}) \subset M^{ur}(\rho_{|I_M}). \end{equation}

\noindent Obviously, these inclusions are equality if and only if $N=0$, that is $\rho=\rho^W$. We used implicitly this trivial fact in the proof of 
lemma \ref{rep} (ii). \ps

Assume now that $l=p$, and that $\rho: \G_M \rightarrow \GL(V)$ is a continuous representation on a finite dimensional $\Qpb$-vector space $V$, 
which is potentially semistable in the sense of Fontaine \cite{F}. It means that 
$\Dpst(V):=\bigcup_{L/K {\rm finite}} (V \otimes_{\Q_p} B_{st})^{\G_L}$
is free of rank $\dim_{\Qpb}(V)$ over $\Qpb \otimes_{\Q_p} \Q_p^{nr}$, where $B_{st}$ is the usual Fontaine ring. Fontaine defines then a representation, say $\rho^W$, of
$\W_M$ on $\Dpst(V)$, whose restriction to $I_M$ is the obvious one. However, 
the analogue of the inclusion (\ref{inclusion}) is not always satisfied in this setting\footnote{It is however obviously satisfied when 
$\rho_{|I_M}$ has finite image, as 
$\rho^W=\rho_{|W_M}$ in this case.}. We are grateful to P. Colmez and T. Saito for explaining us the following example: \ps
\newcommand{\rec}{\mathrm{rec}}
{\it Example:} Let $M:=\Q_p$ and $K \subset \Qpb$ be a finite extension of $\Q_p$. Let $\omega: \G_K \rightarrow \Qpb^*$ 
be a Lubin-Tate character of $K$, and $\psi: \G_K \rightarrow \Qpb^*$ be a
continuous character of finite order on $I_K$, and $\chi:=\omega\psi$. It is well known that $\omega$ is crystalline,
%of weights $(1,0,\dots,0)$ ($1$ appears once, for the canonical inclusion $K \rightarrow \Qpb$), 
hence $\chi$ is $pst$, and the representation
${\chi^W}_{|I_K}$ on $\Dpst(\chi)$ is $\psi_{|I_K}$. Let $\rho:=\Ind_{\G_K}^{\G_{\Q_p}} \chi$.
An easy computation shows that $\Dpst(\rho)=\Ind_{\W_K}^{\W_{\Q_p}} \Dpst(\chi)$, hence $\rho$ is
$pst$ and $\rho^W=\Ind_{\W_K}^{\W_{\Q_p}}\chi^W$.
\par
        Let us assume from now that $K/\Q_p$ is Galois to simplify. Let $\rec: \OO_K^* \rightarrow \G_K^{ab}$ be 
the reciprocity map of local classfield theory, recall that $\omega \rond \rec$ is
the natural map $\mathcal{O}_K^* \rightarrow \Qpb^*$ induced by some
field embedding $K \rightarrow \Qpb$. We choose now $\psi$ such that $\psi
\rond \rec$ is trivial on the pro-$p$-Sylow of $\OO_K^*$, and coincides with 
$(\omega \rond \rec)^{-1}$ on its $p'$-part. Then, $F_1:=\Q_p^{ur}(\rho^W)$ is the
maximal, tamely ramified, abelian extension of $K$, and is linearly
disjoint with $F_2:=\Q_p^{ur}(\rho)$ over $K^{ur}$. In particular, when
$K/\Q_p$ is unramified of degree $>1$, it gives examples of 
$p$-adic $pst$ representations $\rho$ of $\G_{\Q_p}$ such that $\Q_p^{ab}(\rho^W)$ is not included in $\Q_p^{ab}(\rho)$. \medskip

{\it Remarks: } (i) We do not know if there exists a strict subfield $L$ of
$\Qpb$ such that the inclusion (\ref{inclusion}) holds over $L$ for any $\rho$, or if (\ref{inclusion})
(or some variant) holds under a mild hypothesis on $\rho$. For our aim,
the only controls we seem to have on $\rho$ are $\rho^W_{|I_M}$ and some information on its set of 
Hodge-Tate weights. What we do not control at all is the filtration data of
$\Dpst(\rho)$. In this direction, we may hope that the $p$-adic Langlands philosophy initiated by Breuil 
will allow, in the future, to use the Galois representations attached to $p$-adic automorphic forms, 
whose properties at $p$ are less restricted. \ps
	(ii) Assuming the extension of T. Saito's theorem \cite{S} to
Harris-Taylor's automorphic forms, and using corollary \ref{corgeneral}, a
solution of the above problem of kernels would imply property
$P_{E,\{\infty,u,u'\},u}$ where $E/\Q$ is a quadratic imaginary field and
$p=uu'$ a prime number which splits in $E$. We mention this as the hypothesis we
made in \S \ref{conj} still seem to be out of reach (especially hypothesis 2).

\bigskip
\medskip

\begin{center}{\sc acknowledgments:} \end{center}

The author thanks J. Bella\"iche, Y. Benoist, L. Clozel, P. Colmez, J.-F. Dat, L. Fargues, A. Genestier, M. Harris, G. Henniart, J.-P. Labesse,
T. Saito, and K. Wingberg, for their remarks or helpful discussions. \ps

\bigskip

\noindent CHENEVIER Ga\"etan, \par \noindent
CNRS U.M.R. 7539, L.A.G.A., Universit\'e Paris 13, \pn 
99 Av. J-B. Cl\'ement, 93430 Villetaneuse, France, \par
\noindent E-mail: gaetan.chenevier@normalesup.org

\end{document}